\documentclass{amsart}
\usepackage{amsxtra}
\usepackage{mathtools}
\usepackage{mathdots}
\usepackage{tensor}
\usepackage{multirow}
\usepackage{tikz-cd}
\usepackage{ytableau}
\usepackage{mathrsfs}
\usepackage{verbatim}
\usepackage{enumitem}
\usepackage{graphicx}
\usepackage{calc}
\usepackage{hyperref}
\usepackage[shortalphabetic,initials]{amsrefs}
\usepackage{bdsstandard}
\usepackage{bdsthmart}

\numberwithin{equation}{section}

\newcommand{\aform}{\ensuremath{\la\,\cdot\, , \,\cdot\,\ra}\xspace}

\newcommand{\DOrb}{\operatorname{\partial Orb}}
\newcommand{\GS}{\ensuremath{\mathrm{GS}}\xspace}
\newcommand{\loc}{\ensuremath{\mathrm{loc}}\xspace}

\newcommand{\sform}{\ensuremath{(\,\cdot\, , \,\cdot\,)}\xspace}

\newcommand{\spl}{\ensuremath{\mathrm{spl}}\xspace}

\renewcommand{\L}{\ensuremath{\mathscr{L}}\xspace}
\renewcommand{\H}{\ensuremath{\mathscr{H}}\xspace}

\newcommand\longvdots[1]{\rotatebox{-90}{\hbox to #1 {\dotfill}}}

\newcommand{\squigra}{%
   {\mathrel{\xymatrix@C=6.4ex@M=.12ex{\ar@{~>}[r] &}}}%
   }
\newcommand{\wtHG}{\ensuremath{{\widetilde{{H\! G}}}}\xspace}

\begin{document}
   
\renewcommand{\O}{\ensuremath{\mathscr{O}}\xspace}

\title{An arithmetic intersection conjecture}
\author{Brian Smithling}
\address{Johns Hopkins University, Department of Mathematics, 3400 N.\ Charles St.,\ Baltimore, MD  21218, USA}
\email{bds@math.jhu.edu}
\thanks{This work is supported by NSA grant H98230-16-1-0024.}

\begin{abstract}
We survey recent joint work with M. Rapoport and W. Zhang related to the arithmetic Gan--Gross--Prasad conjecture for Shimura varieties attached to unitary groups.
\end{abstract}

\maketitle

\section{Introduction}\label{intro}
This article is a survey centered on some recent joint work \cite{rsz3} of the author with M. Rapoport and W. Zhang. It is based on the author's talk at the 2017 First Annual Meeting of the International Consortium of Chinese Mathematicians.

Let us begin by discussing a little of the historical context.
One of the pioneering results in arithmetic geometry in the past 35 years is the theorem of Gross and Zagier \cite{grosszagier86}, which relates the N\'eron--Tate height of Heegner points on modular curves to the first central derivative of certain $L$-series attached to modular forms.  
Many mathematicians have worked to extend and generalize this result in various ways during the intervening years.
One such conjectural generalization to higher dimensional varieties, the \emph{arithmetic Gan--Gross--Prasad (AGGP) conjecture}, was proposed by Gan, Gross, and Prasad in \cite{ggp12a}, in the context of certain Shimura varieties attached to unitary and orthogonal groups.  Now, around the same time as \cite{grosszagier86} appeared, Waldspurger \cite{wald85} obtained a similar formula for the central value of the $L$-series in question, in terms of toric period integrals. Waldspurger's formula was subsequently conjecturally generalized to special orthogonal groups by Gross and Prasad \cite{gropra92,gropra94}, and to broader classes of classical groups by Gan, Gross, and Prasad \cite{ggp12a}.  In \cite{jacral11}, Jacquet and Rallis proposed a \emph{relative trace formula} (RTF) approach to the resulting Gan--Gross--Prasad conjecture in the case of unitary groups.  In \cite{zhang12a}, Zhang proposed a kind of ``differentiated'' version of this approach in the context of the unitary case of the AGGP conjecture.  The comparison entailed in his approach gave rise to a local conjecture, the \emph{arithmetic fundamental lemma} (AFL), which has inspired much activity in recent years. For example, we mention direct work on the AFL by Rapoport--Terstiege--Zhang \cite{rtz13}, Mihatsch \cite{mihatsch17,mihatschTh}, and Li--Zhu \cite{lizhu17}; work on some variants where ramification is allowed in \cite{rsz17,rsz18}; and work on a variant for Lubin--Tate spaces by Li \cite{li?}.  By contrast, the paper \cite{rsz3} returns to the global setting of the AGGP conjecture as developed in \cite{zhang12a}.

In fact, Zhang's approach led him to 
develop
 a variant of the AGGP conjecture which relates intersection numbers of cycles to derivatives of the distributions appearing in the Jacquet--Rallis RTF.  One of the main goals of \cite{rsz3} is to extend and complete this variant to a full global conjecture whose proof should be accessible (although still a challenge!)\ in various cases.  Indeed, even in its formulation, the intersection-theoretic side of the AGGP conjecture relies on the existence of the Beilinson--Bloch height pairing, which is contingent on some widely open conjectures on algebraic cycles.  By contrast, the arithmetic intersection conjecture in \cite{rsz3} can be stated unconditionally (it defines intersection numbers in terms of the Gillet--Soul\'e pairing on arithmetic Chow groups). We remark that the distributions in the Jacquet--Rallis RTF are related to the $L$-functions that appear in the AGGP conjecture, which justifies that the conjecture of \cite{rsz3} can be considered as a variant of the AGGP conjecture.

Another main aim of \cite{rsz3} is to shift the setting of the conjecture to variants of the Gan--Gross--Prasad Shimura varieties which are of PEL type, i.e., which admit moduli descriptions in terms of abelian schemes equipped with a polarization, endomorphisms, and a level structure.  The resulting Shimura varieties give finite \'etale covers of the Gan--Gross--Prasad Shimura varieties, and they carry the advantage of admitting natural integral models.

The aim of this survey is to serve as an introduction to the above aspects of \cite{rsz3}, focusing most of the time on the case of unitary groups attached to an imaginary quadratic field.  Many aspects of this material are also surveyed by Zhang in \cite{zhang12b,zhangicm}.  Here is an outline of the contents.  In \s\ref{aggp}, we recall (a version of) the arithmetic Gan--Gross--Prasad conjecture.  In \s\ref{rsz} we introduce the Shimura varieties of \cite{rsz3} and their so-called semi-global integral models (i.e.,\ models defined over the localization of the ring of integers of the reflex field at a prime).  In \s\ref{sglob}, we formulate the semi-global version of the arithmetic intersection conjecture for these integral models; this is both simpler than the full global version, and it is also the version for which we have concrete evidence.  We conclude the paper in \s\ref{glob} by briefly discussing the integral models of our Shimura varieties and the arithmetic intersection conjecture in the global case.

Throughout the paper, we write $\AA_F$, $\AA_{F,f}$, and $\AA_{F,f}^p$ for the respective rings of adeles, finite adeles, and finite adeles away from $p$ of a number field $F$. We abbreviate these to \AA, $\AA_f$, and $\AA_f^p$ when $F = \QQ$.

\subsection*{Acknowledgments} 
It is a pleasure to thank the ICCM for the opportunity to present this material at the 2017 meeting.  I also thank M. Rapoport, W. Zhang, and T. Haines for helpful comments.

\section{The arithmetic Gan--Gross--Prasad conjecture}\label{aggp}
We begin in this section by stating a version of the unitary case of the arithmetic Gan--Gross--Prasad conjecture \cite{ggp12a}.  

\subsection{GGP Shimura varieties}\label{ggp sh var}
Although later in the paper we will restrict to the case of an imaginary quadratic field in the setup, for now it causes little harm to allow an arbitrary CM field.  So let $F$ be a CM number field, that is, a totally imaginary quadratic extension of a totally real number field $F_0$.  Let $d := [F_0 : \QQ]$.  Let $\Phi \subset \Hom_\QQ(F,\CC)$ be a CM type for $F$, that is, a set of $d$ embeddings whose restrictions to $F_0$ yield the $d$ distinct embeddings $F_0 \to \RR$.  Let $n$ be a positive integer, and let $(W,h)$ be an $n$-dimensional $F/F_0$-hermitian space whose signatures at the archimedean places are of the form
\[
   \sig (W \otimes_{F_0,\varphi} \RR) = 
   \begin{cases} 
      (1,n-1),  &  \text{$\varphi = \varphi_0$ for some distinguished $\varphi_0 \in \Phi$};\\
      (0,n),  &  \varphi \in \Phi \ssm \{\varphi_0\}.
   \end{cases}
\]
Let $\rU(W)$ denote the unitary group of $W$ (a reductive group over $F_0$), and let
\[
   G := \Res_{F_0/\QQ}\rU(W).
\]
In other words, $G$ is the unitary group of $W$ regarded as an algebraic group over \QQ; in particular, $G(\QQ) = \rU(W)(F_0)$.

Attached to the above data is a natural Shimura datum $(G,X)$, where $X$ is a $G(\RR)$-conjugacy class of cocharacters $\Res_{\CC/\RR}\GG_m \to G_\RR$ which may be identified as a Hermitian symmetric domain with the $(n-1)$-dimensional open unit ball $\{ z \in \CC^{n-1} \mid z \cdot \ol z < 1\}$; see \cite{ggp12a}*{\s27} or \cite{rsz3}*{Rems.\ 3.2(iii), 3.3}.  By the formalism of Shimura varieties, we obtain an inverse system $(\Sh_K(G,X))_K$ of smooth quasi-projective varieties over \CC indexed by sufficiently small compact open subgroups $K \subset G(\AA_f)$.  In fact, for arbitrary compact open $K$ we may still regard $\Sh_K(G,X)$ in a natural way as a smooth Deligne--Mumford stack,\footnote{Although, abusing language, we will often still refer to such objects as Shimura varieties.} and in all cases there is an identification of complex points
\[
   \Sh_K(G,X)(\CC) \cong G(\QQ) \big\bs \bigl[X \times G(\AA_f) / K\bigr],
\]
where the right-hand side is interpreted as the groupoid quotient.  Again by the general formalism, these Shimura varieties admit \emph{canonical models} over a number field in \CC prescribed by the Shimura datum (the \emph{reflex field}), which in the case at hand identifies with $F$ embedded via $\varphi_0$.

When $F_0 \neq \QQ$, the existence of a $\varphi$ at which $W \otimes_{F_0,\varphi} \RR$ is a definite Hermitian space implies that $G$ is anisotropic over \QQ, and hence $\Sh_K(G,X)$ is compact.  When $F_0 = \QQ$, these Shimura varieties are typically not compact, but each possesses a unique toroidal compactification \cite{amrt10,pink90}.  
Abusing notation, we will understand $\Sh_K(G,X)$ to 
designate
the replacement of the Shimura variety by this compactification whenever the former is non-compact.

\subsection{Beilinson--Bloch height pairing}
Let $Y$ be a smooth proper variety of dimension $d$ over a number field $E$.  For each $0 \leq i \leq d$, let $\Ch^i(Y)$ denote the Chow group on $Y$ of codimension-$i$ algebraic cycles modulo rational equivalence.  We consider the cycle class map valued in Betti cohomology,
\[
   \cl_i \colon \Ch^i(Y)_\QQ \to \rH^{2i}\bigl(Y(\CC),\QQ\bigr).
\]
Its kernel is the space of \emph{cohomologically trivial cycles} in degree $i$,
\[
   \Ch^i(Y)_{\QQ,0} := \ker \cl_i.
\]
Contingent on some widely open conjectures on algebraic cycles and the existence of regular proper integral models of $Y$, Beilinson \cite{bei87} and Bloch \cite{bloch84} have defined a \emph{height pairing}
\[
   \sform_{\mathrm{BB}}\colon \Ch^i(Y)_{\QQ,0} \times \Ch^{d-i+1}(Y)_{\QQ,0} \to \RR.
\]
When $i = 1$, the height pairing is unconditionally defined and coincides with the classical N\'eron--Tate height pairing between divisors and zero-cycles.  The height pairing is also defined whenever $Y$ admits a smooth proper model over $\Spec O_E$ \cite{rsz3}*{Rem.\ 6.4} (even in the case that $Y$ and its integral model are DM stacks).

\subsection{AGGP Conjecture}\label{ss:aggp conj}
We return to the setting of \s\ref{ggp sh var} and now assume that $n \geq 2$. Fix a totally negative vector $u \in W$, i.e., a vector such that $\varphi(h(u,u)) < 0$ for all $\varphi \in \Phi$.  Let $W^\flat := (u)^\perp \subset W$.  The restriction of $h$ to $W^\flat$ makes $W^\flat$ into a Hermitian space of the same signature type at the archimedean places as $W$, with $n-1$ in place of $n$.  Therefore we may apply the discussion in \s\ref{ggp sh var} to $W^\flat$: we define
\[
   H := \Res_{F_0/\QQ}\rU(W^\flat),
\]
and we obtain a natural Shimura datum $(H, X^\flat)$.  The group $H$ embeds canonically in $G$ as the stabilizer of $u$, and this embedding induces a morphism of Shimura data $(H, X^\flat) \to (G, X)$.  Therefore, whenever there is an inclusion of compact open subgroups $K_H \subset K_G$, we obtain a (finite and unramified) morphism of Shimura varieties,
\begin{equation}\label{morph}
   \Sh_{K_H}(H) \to \Sh_{K_G}(G)
\end{equation}
(here and below we suppress $X^\flat$ and $X$ in the notation).  We then consider the graph morphism
\begin{equation}\label{graph}
   \Delta_{K_H,K_G}\colon \Sh_{K_H}(H) \inj \Sh_{K_H}(H) \times \Sh_{K_G}(G) = \Sh_{K_H \times K_G}(H \times G).
\end{equation}
Here the source has dimension $n-2$ and the target has dimension $2n-3$, so that we obtain a cycle ``just below'' the middle dimension of the target.  The morphism of Shimura data $(H,X^\flat) \to (H \times G, X^\flat \times X)$ is therefore a \emph{special pair} in the sense of Zhang \cite{zhangicm}*{Def.\ 3.1}.  We obtain a class in the Chow group,
\[
   z_{K_H,K_G} := (\Delta_{K_H,K_G})_*[\Sh_{K_H}(H)] \in \Ch^{n-1}\bigl(\Sh_{K_H \times K_G}(H \times G)\bigr)_\QQ.
\]

Now, the Beilinson--Bloch height pairing on a $(2n-3)$-dimensional variety gives a self-pairing on the codimension-$(n-1)$ cohomologically trivial Chow group, and, roughly speaking, the intersection-theoretic side of the AGGP conjecture is in terms of this pairing against the cycle class $z_{K_H,K_G}$.  Of course, for this to honestly make sense, it is necessary to replace $z_{K_H,K_G}$ by a ``cohomological trivialization'' in $\Ch^{n-1}(\Sh_{K_H \times K_G}(H \times G))_{\QQ,0}$.  In \cite{rsz3}*{\s6.3}, this is done by defining a map
\begin{equation}\label{R(f-)}
   \Ch^{n-1}\bigl(\Sh_{K_H \times K_G}(H \times G)\bigr)_\QQ \to \Ch^{n-1}\bigl(\Sh_{K_H \times K_G}(H \times G)\bigr)_{\QQ,0}
\end{equation}
via the action of a certain Hecke correspondence.  We denote the image of $z_{K_H,K_G}$ under \eqref{R(f-)} by $z_{K_H,K_G,0}$, and we call it the \emph{cohomologically trivial arithmetic diagonal cycle}.  Subject to standard conjectures on algebraic cycles, the map \eqref{R(f-)} is a retraction onto the target, and formation of $z_{K_H,K_G,0}$ is compatible with push-forward under change of level subgroup \cite{rsz3}*{Rem.\ 6.9}.  

Now define (with respect to the natural pullback maps)
\[
   \Ch^{n-1}\bigl(\Sh(H \times G)\bigr)_0 := \clim_{\substack{K_H \subset H(\AA_f)\\ K_G \subset G(\AA_f)}} \Ch^{n-1}\bigl(\Sh_{K_H \times K_G}(H \times G)\bigr)_{\CC,0}.
\]
Then the \CC-linear extension of the pairing $\sform_{\mathrm{BB}}$ against the compatible family $(z_{K_H,K_G,0})_{K_H \subset K_G}$ defines a linear functional
\[
   \ell \colon \Ch^{n-1}\bigl(\Sh(H \times G)\bigr)_0 \to \CC.
\]
We arrive at the following version of the arithmetic Gan--Gross--Prasad conjecture \cite{ggp12a}*{\s27}.
\begin{conj}[Arithmetic Gan--Gross--Prasad]\label{aggp conj}
Let $\pi$ be a tempered cuspidal automorphic representation of $(H \times G)(\AA)$ appearing in the cohomology $\rH^*(\Sh(H \times G))$. The following are equivalent.
\begin{enumerate}[label=(\alph*),ref=\alph*]
\item The restriction of $\ell$ to the $\pi_f$-isotypic component of $\Ch^{n-1}(\Sh(H \times G))_0$ is nonzero.
\item\label{ggp ii} $\Hom_{H(\AA_f)}(\pi_f,\CC) \neq 0$ and the first order central derivative $L'(\frac 1 2, \pi, R) \neq 0$.
\end{enumerate}
\end{conj}

Here $\rH^*(\Sh(H \times G)) := \clim_K \rH^*(\Sh_K(H \times G)(\CC),\CC)$, and the $L$-function appearing in \eqref{ggp ii} is the Rankin--Selberg convolution $L$-function of the base change of $\pi$ to $(\GL_{n-1} \times \GL_n)(\AA_F)$ (see e.g.\ \cite{rsz3}*{\s6.4}).  In the case $n = 2$ (i.e.,\ the case that $\Sh(H \times G)$ is a curve), the Beilinson--Bloch height pairing is well--defined and coincides with the N\'eron--Tate pairing, and Conjecture \ref{aggp conj} is proved by Yuan--Zhang--Zhang \cite{yzz13}.  But we know of no higher-dimensional cases in which the conjecture is proved.\footnote{However, the \emph{orthogonal group} version of the conjecture has been proved by Yuan--Zhang--Zhang \cite{yzz?} in some special cases where the ambient Shimura variety has dimension $3$.}

\section{RSZ Shimura varieties}\label{rsz}
In this section we introduce the Shimura varieties of \cite{rsz3}, which are variants of those of Gan--Gross--Prasad.  We continue with the notation of \s\ref{aggp}, but to simplify a number of technical aspects, we restrict to the case $F_0 = \QQ$; thus $F$ is imaginary quadratic and $\Phi = \{\varphi_0\}$.  See \s\ref{F_0 neq Q} for remarks on the case of more general $F_0$.  We denote the nontrivial automorphism of $F$ by $a \mapsto \ol a$.

\subsection{Shimura data}\label{rsz sh data}
Let $Z$ be the \QQ-torus
\[
   Z := \Res_{F/\QQ}\GG_m.
\]
Using $\varphi_0$, we identify $\smash[b]{F \otimes \RR \xra[\undertilde]{\varphi_0 \otimes 1} \CC}$ and hence $Z_\RR \isoarrow \Res_{\CC/\RR}\GG_m$, and we define the Shimura homomorphism
\[
   h_{Z}\colon
   \begin{tikzcd}[baseline=(x.base),row sep=-.5ex]
      |[alias=x]| \Res_{\CC/\RR}\GG_m \ar[r]  &  Z_\RR\\
      z \ar[r, mapsto]  &  \ol z.
   \end{tikzcd}
\]
Then $(Z, \{h_Z\})$ is a Shimura datum with reflex field $F$, over which the canonical models of the resulting Shimura varieties are finite and \'etale.

The variants of the Gan--Gross--Prasad Shimura data we will consider are now obtained by everywhere taking the product with $(Z,\{h_Z\})$: we define the groups
\[
   \wt G := Z \times G,
   \quad
   \wt H := Z \times H,
   \quad
   \wtHG := Z \times H \times G,
\]
and the Shimura data
\[
   \bigl(\wt G, \{h_Z\} \times X\bigr),
   \quad
   \bigl(\wt H, \{h_Z\} \times X^\flat\bigr),
   \quad
   \bigl(\wtHG, \{h_Z\} \times X^\flat \times X\bigr).
\]
Thus we obtain product decompositions of pro-varieties,
\begin{gather*}
   \Sh(\wt G) = \Sh(Z) \times \Sh(G),
   \quad
   \Sh(\wt H) = \Sh(Z) \times \Sh(H),\\
   \Sh(\wtHG) = \Sh(Z) \times \Sh(H) \times \Sh(G)
\end{gather*}
(as before we suppress 
$\{h_Z\}$, $X^\flat$, and $X$ in the notation).  In particular, the Shimura varieties for the ``tilde'' groups give finite \'etale covers of the Gan--Gross--Prasad Shimura varieties.  All of these Shimura varieties have common reflex field $F$.

\subsection{Semi-global integral models}\label{sglob int mod}
We are now going to describe ``semi-global'' integral models (i.e.,\ models over the localization of $O_F$ at a prime ideal) of the Shimura varieties of \s\ref{rsz sh data} for certain level subgroups, in terms of explicit moduli problems.

In the case of the group $Z$, the Shimura variety has a well-known \emph{global} integral model over $\Spec O_F$, which we will recall only in the case of the (unique) maximal compact subgroup
\[
   K_Z^\circ := \wh O_F^\times \subset Z(\AA_f) = \AA_{F,f}^\times.
\]
For each $O_F$-scheme $S$, let $\M_0(S)$ be the groupoid of pairs $(A_0,\iota_0)$, where $A_0$ is an elliptic curve over $S$ and $\iota_0\colon O_F \to \End_S(A_0)$ is an action whose induced action on $\Lie A_0$ is the composition of the nontrivial Galois automorphism of $O_F$ with the natural action arising from the structure map $O_F \to \O_S$.  An isomorphism in $\M_0(S)$ is an $O_F$-linear isomorphism of elliptic curves over $S$.  Then $\M_0$ is representable by a Deligne--Mumford stack finite and \'etale over $\Spec O_F$, and its coarse moduli scheme is isomorphic to $\Spec O_H$, where $H$ denotes the Hilbert class field of $F$ \cite{kry99}*{\s5}.  The generic fiber of $\M_0$ is the canonical model of $\Sh_{K_Z^\circ}(Z)$ \cite{kudrap14a}*{Props.\ 4.3, 4.4}.

To define semi-global models for $\Sh(\wt G)$, we first specify integral data.  Fix a prime number $p$, and assume that $p$ splits in $F$ if $p = 2$.  Set $F_p := F \otimes_\QQ \QQ_p$ and $W_p := W \otimes_\QQ \QQ_p$ (an $F_p/\QQ_p$-hermitian space).  Let $\pi_p$ be a uniformizer in $F_p$ (we simply take $\pi_p = p$ if $p$ is inert or split in $F$).   We choose an $O_{F_p}$-lattice $\Lambda_p \subset W_p$ of the following form:
\begin{itemize}
\item (self-dual)  $\Lambda_p = \Lambda_p^\vee$  if $p$ splits in $F$, or if $p$ is inert in $F$ and $W_p$ is a split Hermitian space;
\item (almost self-dual)  $\Lambda_p \subset^1 \Lambda_p^\vee \subset p\i\Lambda_p$ if $p$ is inert in $F$ and $W_p$ is a non-split Hermitian space;
\item ($\pi_p$-modular)  $\Lambda_p^\vee = \pi_p\i\Lambda_p$ if $p$ ramifies in $F$ and $n$ is even;
\item (almost $\pi_p$-modular)  $\Lambda_p \subset \Lambda_p^\vee \subset^1 \pi_p\i\Lambda_p$ if $p$ ramifies in $F$ and $n$ odd.
\end{itemize}
Here $\Lambda_p^\vee = \{x \in W_p \mid h(x,\Lambda_p) \subset O_{F_p}\}$ denotes the dual lattice, and the symbol $\subset^1$ means that the inclusion is of $O_{F_p}$-colength one.  When $p$ is inert in $F$, there are two $n$-dimensional $F_p/\QQ_p$-Hermitian spaces up to isometry (\emph{split} and \emph{non-split}), and they are distinguished by whether they respectively contain a self-dual or an almost self-dual lattice.  When $p$ ramifies in $F$, there are again two $n$-dimensional $F_p/\QQ_p$-Hermitian spaces up to isometry; when $n$ is odd, both contain an almost $\pi_p$-modular lattice, but when $n$ is even, the existence of a $\pi_p$-modular lattice determines the isometry type of $W_p$.  Thus we are imposing a condition on $W_p$ in the third bullet.  We consider a level subgroup of the form
\[
   K_{\wt G} := K_Z^\circ \times \Stab(\Lambda_p) \times K_G^p \subset \wt G(\AA_f) = Z(\AA_f) \times G(\QQ_p) \times G(\AA_f^p),
\]
where $K_G^p \subset G(\AA_f^p)$ is arbitrary.  We note that $\Stab(\Lambda_p)$ is a maximal parahoric subgroup of $G(\QQ_p)$ in the first three bullets above, and it contains a maximal parahoric subgroup with index two in the fourth bullet.

Now let $w$ be a place of $F$ over $p$.  For each locally Noetherian $O_{F,(w)}$-scheme $S$, let $\M_{K_{\wt G}}(\wt G)(S)$ be the groupoid of tuples $(A_0,\iota_0,A,\iota,\lambda,\ol\eta^p)$, where
\begin{itemize}
\item $(A_0,\iota_0) \in \M_0(S)$;
\item $(A,\iota)$ is an abelian scheme over $S$ with an action $O_F\otimes\ZZ_{(p)} \to \End_S(A) \otimes \ZZ_{(p)}$ satisfying the \emph{Kottwitz condition} 
\begin{equation}\label{kott}
   \charac\bigl(\iota(a) \mid \Lie A \bigr) = (T-a)(T- \ol a)^{n-1}, \quad a \in O_F \otimes \ZZ_{(p)};
\end{equation}
\item $\lambda\colon A \to A^\vee$ is a polarization whose Rosati involution induces the nontrivial Galois automorphism on $O_F \otimes \ZZ_{(p)}$, and such that $\ker\lambda \subset A[\pi_p]$ of rank $\#(\Lambda_p^\vee/\Lambda_p)$; and
\item $\ol\eta^p$ is a $K_G^p$-orbit of isometries of $\AA_{F,f}^p/\AA_{F_0,f}^p$-Hermitian spaces,
\begin{equation}\label{level}
   \eta^p\colon \Hom_{\AA_{F,f}^p}\bigl(\rV^p(A_0),\rV^p(A)\bigr) \isoarrow (-W) \otimes_F \AA_{F,f}^p.
\end{equation}
\end{itemize}
Here $-W$ denotes the Hermitian space $(W,-h)$, the notation $\rV^p(\,\cdot\,)$ denotes the prime-to-$p$ rational Tate module, and the $\Hom$ space in \eqref{level} carries a natural Hermitian form arising from $\lambda$ and the unique principal polarization $\lambda_0$ on $A_0$; see \cite{rsz3} for details.  To complete the moduli problem, when $p$ is ramified, it is also necessary to impose the \emph{wedge condition}
\[
   \bigwedge\nolimits^2 \bigl(\iota(a) - \ol a \mid \Lie A \bigr) = 0, \quad a \in O_F \otimes \ZZ_{(p)},
\]
and the \emph{(refined) spin condition} on $\Lie A$; this last condition is technical to formulate (at least when $n$ is odd) and we refer to \cite{rsz3}*{\s4.4} and the references therein for details.  An isomorphism $(A_0, \iota_0, \lambda_0, A, \iota, \lambda, \ol\eta^p) \isoarrow (A_0', \iota_0', \lambda_0', A', \iota', \lambda', \ol\eta'^p)$ in $\M_{K_{\wt G}}(\wt G)(S)$ consists of an isomorphism $(A_0,\iota_0) \isoarrow (A_0',\iota_0')$ in $\M_0(S)$ and a $p$-principal $O_F \otimes \ZZ_{(p)}$-linear quasi-isogeny $A \to A'$ such that $\lambda'$ pulls back to $\lambda$ and $\ol\eta'^p$ pulls back to $\ol\eta^p$.  We have the following structure result from \cite{rsz3}.

\begin{thm}\label{semi-glob model thm}
$\M_{K_{\wt G}}(\wt G)$ is a Deligne--Mumford stack of relative dimension $n-1$ over $\Spec O_{F,(w)}$, smooth except in the ``almost self-dual'' case, where it is of semi-stable reduction.  Furthermore, the generic fiber $\M_{K_{\wt G}}(\wt G)_F$ is naturally isomorphic to the canonical model of $\Sh_{K_{\wt G}}(\wt G)$.
\end{thm}

In particular, one has the surprising fact that $\M_{K_{\wt G}}(\wt G)$ is smooth even when $p$ ramifies in $F$, a phenomenon termed ``exotic smoothness'' in \cite{rsz17}.  The proofs of this and all other statements on the local structure of $\M_{K_{\wt G}}(\wt G)$ in Theorem \ref{semi-glob model thm} reduce to analogous statements for the \emph{local model} $M^\loc$ \cite{rapzink96,prs13}, which are due to various authors in various cases.  Roughly speaking, the local model is a moduli space for linear-algebraic data modeling the Hodge filtration $\rH_1^{\mathrm{dR}}(A) \surj \Lie A$ of the abelian variety in the moduli problem.  For example, in the ``self-dual'' case, $M^\loc$ boils down to the moduli problem of line bundle quotients $\O^n \surj \L$; thus $M^\loc \cong \PP^{n-1}$, and the local assertions visibly hold.  We note that He--Pappas--Rapoport \cite{hpr?} have recently made substantial progress towards classifying the cases in which a general local model has good or semi-stable reduction.

\begin{rem}\label{drin}
The above Shimura variety for $\wt G$ and its moduli interpretation (in the case $F_0 = \QQ$) is closely related to the one studied in \cite{bhkry?} by Bruinier--Howard--Kudla--Rapoport--Yang.  In fact, loc.\ cit.\ defines a global integral model over $\Spec O_F$ for (essentially) the same Shimura datum, under the assumption that $W$ contains a self-dual lattice $\Lambda$ and the level subgroup $K_{\wt G}$ is $K_Z^\circ \times \Stab(\Lambda) \subset \wt G(\AA_f) = Z(\AA_f) \times G(\AA_f)$.  If $p$ is split or inert, then the base change of their model to $\Spec O_{F,(w)}$ is of the form of ours above.  But their model is different above ramified primes $p$; in particular, it is not regular above ramified $p$ when $n > 2$.  If $2$ splits in $F$, and if $W$ is split at every ramified prime in the case that $n$ is even, then the semi-global models we have defined can similarly be glued to form a global model over $\Spec O_F$ which has everywhere good or semi-stable reduction (cf.\ \s\ref{glob}), and it would be interesting to consider this stack in the context of \cite{bhkry?}.
\end{rem}

\begin{rem}
In the case that $p$ splits in $F$, one can can add a \emph{Drinfeld level structure} at $w$ to the above moduli problem to obtain a semi-global integral model for the level subgroup where the $p$-factor $\Stab(\Lambda_p)$ is replaced by a principal congruence subgroup in $G(\QQ_p)$ for $\Lambda_p$.  The resulting stack is regular \cite{rsz3}*{Th.\ 4.7}.  This has important technical applications in the context of formulating arithmetic intersection conjectures.
\end{rem}

Of course, there is an analogous semi-global model $\M_{K_{\wt H}}(\wt H)$ for $\Sh_{K_{\wt H}}(\wt H)$, relative to the choice of an appropriate lattice $\Lambda_p^\flat \subset W_p^\flat$.  Setting
\begin{equation}\label{subgroups}
\begin{aligned}
   K_H &:= \Stab(\Lambda_p^\flat) \times K_{\wt H}^p \subset H(\AA_f) = H(\QQ_f) \times H(\AA_f^p),\\
   K_G &:= \Stab(\Lambda_p) \times K_{\wt G}^p \subset G(\AA_f) = G(\QQ_p) \times G(\AA_f^p),\\
   K_{\wtHG} &:= K_Z^\circ \times K_H \times K_G \subset \wtHG(\AA_f) = Z(\AA_f) \times H(\AA_f) \times G(\AA_f),
\end{aligned}
\end{equation}
we also get a model $\M_{K_{\wtHG}}(\wtHG) := \M_{K_{\wt H}}(\wt H) \times_{\M_0} \M_{K_{\wt G}}(\wt G)$ of $\Sh_{K_{\wtHG}}(\wtHG)$.

\subsection{Remarks on the case $F_0 \neq \QQ$}\label{F_0 neq Q}
In this subsection we give a brief guide to the modifications to the above discussion that are required when $F_0 \neq \QQ$.  
This material is somewhat more technical in nature and intended more for experts; 
the reader interested only in the case $F_0 = \QQ$ (which we will continue to assume in later sections) may safely skip over it.

For general $F_0$, one generalizes $Z$ to the group
\[
   Z^\QQ := \bigl\{\, z \in \Res_{F/\QQ} \GG_m \bigm| \Nm_{F/F_0}(z) \in \GG_m \,\bigr\}.
\]
The CM type $\Phi$ determines a standard Shimura homomorphism for $Z^\QQ$, and as before, we take $h_{Z^\QQ}$ to be the precomposition of this homomorphism by complex conjugation.  The Shimura datum $(Z^\QQ,\{h_{Z^\QQ}\})$ has reflex field $E_\Phi$, the reflex field of $\Phi$.  One then defines the groups $\wt G$, $\wt H$, and $\wtHG$ and their Shimura data by taking the product of $(Z^\QQ,\{h_{Z^\QQ}\})$ with the Gan--Gross--Prasad Shimura data.  The resulting three Shimura data have a common reflex field $E$, which is the compositum of $F$ (embedded into \CC via $\varphi_0$) with $E_\Phi$.  In particular, $E$ may be strictly larger than $F$, in constrast to the case $F_0 = \QQ$.

One defines a (global) integral model for $\Sh_{K_{Z^\QQ}^\circ}(Z^\QQ)$, where $K_{Z^\QQ}^\circ$ is the unique maximal compact open subgroup of $Z^\QQ(\AA_f)$, similarly to before, in terms of abelian schemes with $O_F$-action satisfying a Kottwitz condition of signature $((0,1)_{\varphi \in \Phi})$.  Two notable subtleties arise in the general case.  One is that a polarization must be specified as part of the data of this moduli problem, but if $F/F_0$ is unramified at all finite places, then the moduli problem may be empty if one requires the polarization to be principal.  As a remedy, following Howard \cite{how12}*{\s3.1}, one fixes a nonzero ideal $\aaa \subset O_{F_0}$ and defines $\M_0^\aaa$ to represent the version of the moduli problem where the kernel of the polarization is required to equal $A_0[\aaa]$.  Then $\M_0^\aaa$ is finite and \'etale over $\Spec O_{E_\Phi}$, and it is nonempty for appropriate \aaa.  The second subtlety is that, due to the possible failure of the Hasse principle for $Z^\QQ$, the generic fiber of $\M_0^\aaa$ may only be a finite disjoint union of copies of $\Sh_{K_{Z^\QQ}^\circ}(Z^\QQ)$ when it is nonempty.

Now let $\nu$ be a finite place of $E$.  Similarly to before, one defines a semi-global model $\M_{K_{\wt G}}(\wt G)$ over $\Spec O_{E,(\nu)}$ as a moduli space of pairs of abelian schemes $(A_0,A)$ with additional structure, where $A_0$ is required to lie on some fixed summand of $\M_0^\aaa$.  
Let $p$ denote the rational prime under $\nu$, and let $\V_p$ denote the set of places of $F_0$ above $p$.  Of course $\V_p$ contains a distinguished place $v_0$ lying under $\nu$ via $\varphi_0$, but the possible presence of other elements $v \in \V_p$ is an important new source of complications in the moduli problem.
Indeed, this already comes up in the basic question of what kind of level subgroup $K_{\wt G}$ to allow.  If we declare that $K_{\wt G}$ be of the form $K_{Z^\QQ}^\circ \times K_G \subset \wt G(\AA_f) = Z^\QQ(\AA_f) \times G(\AA_f)$,
then it is natural to take the factor $K_G$ to be a product of subgroups according to the decomposition
\[
   G(\AA_f) = G(\QQ_p) \times G(\AA_f^p) = \Biggl( \prod_{v \in \V_p} \rU(W)(F_{0,v}) \Biggr) \times G(\AA_f^p).
\]
Away from $p$, one may take an arbitrary factor $K_G^p \subset G(\AA_f^p)$ as before. At the factors indexed by $v \in \V_p$, one possibility is to choose an appropriate lattice $\Lambda_v \subset W_v$ for each $v$, and then take the corresponding stabilizer group in each factor.  But it is also natural to ask to impose a deeper level structure in the factors indexed by $v \neq v_0$.  Since each such $v$ is of characteristic $p$, the $v$-divisible groups $A_0[v^\infty]$ and $A[v^\infty]$ of the universal abelian schemes are not \'etale over the special fiber of the moduli problem, and therefore one cannot apply the usual notion of level structure (i.e.,\ in terms of Tate modules) in a direct way to them.  However, it is observed in \cite{rsz3}*{Rem.\ 4.2} that the rational Hom module from $A_0[v^\infty]$ to $A[v^\infty]$ forms a lisse $F_{0,v}$-local system, and by considering trivializations of this, one can formulate a notion of $K_G^{v_0}$-level structure for arbitrary $K_G^{v_0} \subset \prod_{v \in \V_p \ssm \{v_0\}} \rU(W)(F_{0,v})$.

The existence of places $v \in \V_p \ssm \{v_0\}$ also necessitates some new conditions in the moduli problem: the \emph{Eisenstein condition} on the Lie algebra of $A$ at such places $v$, which is needed for flatness of the moduli space (it is the analog in this context of the condition introduced by Rapoport and Zink in \cite{rapzink17}), and the \emph{sign condition} at such $v$, which ensures that the generic fiber is a single copy of the Shimura variety for $\wt G$ (it is based on a sign invariant which is similar to one introduced by Kudla and Rapoport in \cite{kudrap15}).  See \cite{rsz3}*{\s4.1}.

There are also some ramification assumptions which are imposed in \cite{rsz3}.  If $v_0$ is ramified over $p$, then it is required to split in $F$, say $v_0 = w_0\ol w_0$, with $w_0$ the place in $F$ under $\nu$ via $\varphi_0$; and $A[w_0^\infty]$ is required to be a Lubin--Tate group in the sense of \cite{rapzink17}.  If $v_0$ ramifies in $F$, or if $v_0$ is inert in $F$ and $W_{v_0}$ is non-split, then $v_0$ is required to be of degree one over $p$ (and as in \s\ref{sglob int mod}, $W_{v_0}$ is required to be split if $v_0$ ramifies and $n$ is even).  Relaxing these assumptions (including allowing $W_{v_0}$ to be non-split in the even ramified case) 
is contingent on further progress
on local models.  Under the above assumptions, the semi-global model $\M_{K_{\wt G}}(\wt G)$ typically has good or semi-stable reduction, and we refer to \cite{rsz3}*{\s4} for the precise statements.

\section{Arithmetic intersection conjecture (semi-global case)}\label{sglob}
In this section we formulate the arithmetic intersection conjecture of \cite{rsz3} in the semi-global case.

\subsection{Distributions}\label{distns}
In this subsection we recall the distributions that appear on the analytic side of our arithmetic intersection conjecture.  For this subsection only, we return to the setting of \s\ref{aggp} and allow $F$ to be an arbitrary CM field.

Define the $F_0$-groups
\[
   G' := \Res_{F/F_0}(\GL_{n-1} \times \GL_n), \quad H_1' := \Res_{F/F_0} \GL_{n-1}, \quad H_2' := \GL_{n-1} \times \GL_n.
\]
For a function $f' \in C_c^\infty(G'(\AA_{F_0})$, we define the automorphic kernel function
\[
   K_{f'}(x,y) := \sum_{\gamma \in G'(F_0)} f'(x\i\gamma y), \quad x,y \in G'(\AA_{F_0}).
\]
Then $K_{f'}(x,y)$ is invariant under the action of $G'(F_0)$ in both variables $x$ and $y$.  For $i = 1,2$, fix a Haar measure on $H_i'(\AA_{F_0})$, take the counting measure on $H_i'(F_0)$, and set $[H_i'] := H_i'(F_0)\bs H_i'(\AA_{F_0})$ with the quotient measure.  For $s \in \CC$, the global distribution on $C_c^\infty(G'(\AA_{F_0}))$ is given by
\begin{equation}\label{J}
   J(f',s) := \int_{[H_1']}\int_{[H_2']} K_{f'}(h_1,h_2)\lv \det(h_1) \rv^s \eta(h_2) \, dh_1 \, dh_2,
\end{equation}
where for $h_2 = (x_{n-1},x_n) \in H_2'(\AA_{F_0}) = \GL_{n-1}(\AA_{F_0}) \times \GL_n(\AA_{F_0})$, we set $\eta(h_2) := \eta_{F/F_0}(x_n)$ if $n$ is even and $\eta(h_2) := \eta_{F/F_0}(x_{n-1})$ if $n$ is odd.  Here $\eta_{F/F_0}\colon F_0^\times \bs \AA_{F_0}^\times \to \{\pm 1\}$ is the idele class character associated to the quadratic extension $F/F_0$ by class field theory.  We formally define
\[
   J(f') := J(f',0)
   \quad\text{and}\quad
   \partial J(f') := \frac{d}{ds}\Big|_{s=0} J(f',s).
\]

The integral \eqref{J} may diverge in general, but it converges on a large class of functions $f'$.  To explain this, consider the natural right action $\gamma \cdot (h_1,h_2) = h_1\i\gamma h_2$ of $H_1' \times H_2'$ on $G'$. Define an element $\gamma \in G'(F_0)$ to be \emph{regular} if its stabilizer with respect to this action is of minimal dimension, and \emph{semi-simple} if its orbit with respect to this action is Zariski-closed.  In the case at hand, for an element to be regular semi-simple, it is equivalent that it have Zariski-closed orbit and trivial stabilizer.  We denote by $G'(F_0)_\rs$ the set of regular semi-simple elements, which is Zariski-open and dense.  These notions continue to make sense with $F_{0,v}$ in place of $F_0$ for any place $v$.  We say that a function $f_v' \in C_c^\infty(G'(F_{0,v}))$ \emph{has regular support} if its support is contained in $G'(F_{0,v})_\rs$.  Finally, we say that a pure tensor $f' = \otimes_v f_v' \in C_c^\infty(G'(\AA_{F_0}))$ \emph{has regular support at $v_0$} if $f_{v_0}'$ has regular support.  Then \eqref{J} converges absolutely on any pure tensor $f'$ which is regularly supported at some place \cite{zhang12a}*{Lem.\ 3.2}.

In fact, more is true for such $f'$: denoting by $[G'(F_0)_\rs]$ the set of regular semi-simple orbits in $G'(F_0)$, there is a decomposition
\[
   J(f',s) = \sum_{\ol\gamma \in [G'(F_0)_\rs]} \, \prod_{\substack{\text{places $v$}\\\text{of $F_0$}}} \Orb(\gamma, f_v',s),
\]
where $\gamma$ denotes any representative of the orbit $\ol\gamma$, and the local orbital integral is defined by
\begin{equation}\label{loc orb}
   \Orb(\gamma, f_v', s) := 
	   \int_{(H_1'\times H_2')(F_{0,v})} f_v'(h_1^{-1}\gamma h_2) \lv\det h_1\rv^s \eta(h_2)\, dh_1\, dh_2.
\end{equation}
It follows that there is a decomposition of the derivative
\begin{equation}\label{delJ decomp}
   \partial J(f') = \sum_v \partial J_v(f'),
\end{equation}
where
\[
   \partial J_v(f') := \sum_{\gamma \in [G'(F_0)_\rs]} \DOrb(\gamma,f_v') \cdot \prod_{u \neq v}\Orb(\gamma,f_u',0),
\]
where in turn
\[
   \DOrb(\gamma,f_v') := \frac{d}{ds}\Big|_{s=0} \Orb(\gamma,f_v',s).
\]
It is the distributions $J$ and $\partial J_v$ that will appear on the analytic side of our arithmetic intersection conjecture.

\subsection{Matching and transfer}\label{m&t}
The relation between the distributions on general linear groups of \s\ref{distns} and the unitary groups of our Shimura varieties is facilitated by the group-theoretic notions of orbit matching and function transfer.  For simplicity we return to the assumption that $F_0 = \QQ$.

In analogy with the case of $H_1' \times H_2'$ acting on $G'$ in the previous subsection, consider the right action of $\wt H \times \wt H$ on $\wt G$ by $g \cdot (h_1,h_2) = h_1\i g h_2$.  We again say that an element $g \in \wtHG(\QQ)$ is \emph{regular semi-simple} if its orbit is Zariski-closed and its stabilizer is of minimal dimension with respect to this action; and the same with $\QQ_\ell$ in place of \QQ for any place $\ell$ (even $\ell = \infty$).  Denoting by $\wtHG(\QQ_\ell)_\rs \subset \wtHG(\QQ_\ell)$ the set of regular semi-simple elements, and by $[\wtHG(\QQ_\ell)_\rs]$ the set of their orbits, there is a natural injection
\begin{equation}\label{orbit inj}
   [\wtHG(\QQ_\ell)_\rs] \inj [G'(\QQ_\ell)_\rs];
\end{equation}
see \cite{rsz3}*{\s2.2} and the references therein.  We remark that for any field $k$ of characteristic not $2$ and any quadratic \'etale algebra $A/k$ in place of $F_\ell/\QQ_\ell$, the target in \eqref{orbit inj} decomposes as a disjoint union of images of analogous embeddings over certain equivalence classes of $A/k$-Hermitian spaces; see \cite{zhangicm}*{(5.5)}.  We say that elements $g \in \wtHG(\QQ_\ell)_\rs$ and $\gamma \in G'(\QQ_\ell)_\rs$ \emph{match} if their respective orbits are identified under \eqref{orbit inj}.

The notion of smooth transfer is now formulated in terms of equalities of orbital integrals along matching regular semi-simple orbits.  On the unitary side, for a function $f_\ell \in C_c^\infty(\wtHG(\QQ_\ell))$ and a regular semi-simple element $g\in \wtHG(\QQ_\ell)_\rs$, we define
\[
   \Orb(g,f_\ell) := \int_{(\wt H(\QQ_\ell)\times\wt H(\QQ_\ell)) /\Delta(Z(\QQ_\ell)) } f_\ell(h_1\i g h_2)\, dh_1\, dh_2.
\]
On the general linear side, for a function $f_\ell' \in C_c^\infty(G'(\QQ_\ell))$ and a regular semi-simple element $\gamma\in G'(\QQ_\ell)_\rs$, we define
\[
   \Orb(\gamma,f_\ell') := \Orb(\gamma,f_\ell',0),   
\]
where the right-hand side is the local orbital integral \eqref{loc orb}.  We then say that the functions $f_\ell$ and $f_\ell'$ are \emph{transfers} of each other if for any element $\gamma \in G'(\QQ_\ell)_\rs$, we have
\[
   \omega_\ell(\gamma) \Orb(\gamma,f_\ell') = 
	   \begin{cases}
			\Orb(g,f_\ell), &  \text{whenever $g \in \wtHG(\QQ_\ell)_\rs$ matches $\gamma$};\\
			0,  &  \text{no $g\in \wtHG(\QQ_\ell)_\rs$ matches $\gamma$}.
		\end{cases}
\]
Here $\omega_\ell$ is an explicit \emph{transfer factor} \citelist{\cite{rsz18}*{\s2.4}\cite{zhang14}*{\s2.4}}.  W. Zhang proved in \cite{zhang14}*{Prop.\ 2.5, Th.\ 2.6} that transfers exist for any $\ell < \infty$, in the sense that any $f_\ell \in C_c^\infty(\wtHG(\QQ_\ell))$ admits a transfer $f_\ell'$, and conversely.  In fact, he proved the more general analog where $F_\ell/\QQ_\ell$ may be replaced by any quadratic \'etale algebra over an $\ell$-adic field (still with $\ell<\infty$), and where on the unitary side one may consider tuples of functions indexed by the equivalence classes of Hermitian spaces (of which there are two in this case) that contribute to the disjoint union decomposition of the target in \eqref{orbit inj} mentioned above.  Existence of transfers remains open in the archimedean case, but Xue has recently established a kind of approximate version in \cite{xue?}.

In the global setting, we say that a pure tensor $f=\otimes_\ell f_\ell\in C_c^\infty(\wtHG(\AA_f))$ and a pure tensor $f'=\otimes_\ell f'_
\ell\in C_c^\infty(G'(\AA_f))$ are \emph{transfers} of each other if they are expressible in a way that $f_\ell$ and $f_\ell'$ are transfers of each other for each prime $\ell$.  Zhang's result together with Yun and Gordon's proof of the \emph{fundamental lemma} for $\ell$ sufficiently large (see the discussion after Conjecture \ref{sglob conj} below) immediately implies the existence of transfers in both directions for pure tensors in $C_c^\infty(\wtHG(\AA_f))$ and in $C_c^\infty(G'(\AA_f))$.

Now, the distributions of \s\ref{distns} are defined on the full adelic function space $C_c^\infty(G'(\AA))$. To extend the notion of transfer to this space, we say that pure tensors $f=\otimes_\ell f_\ell\in C_c^\infty(\wtHG(\AA_f))$ and $f'=\otimes_{\ell\leq\infty} f'_
\ell\in C_c^\infty(G'(\AA))$ are \emph{transfers} of each other if the factors $f_\ell'$ of $f'$ are expressible in a way that $f$ and $\otimes_{\ell < \infty} f_\ell'$ are transfers of each other, and $f_\infty'$ is of a distinguished form, namely a \emph{Guassian function} in the sense of \cite{rsz3}*{Def.\ 7.9}.  The existence of Gaussian functions is still conjectural.

\subsection{Semi-global arithmetic diagonal cycle}\label{sglob adc}
We return to the setup of \s\ref{sglob int mod}.  In this subsection we describe the analog of the morphism of Shimura varieties \eqref{morph} and its graph \eqref{graph} for the semi-global integral models of the RSZ Shimura varieties, defined over $\Spec O_{F,(w)}$ for a $p$-adic place $w$.

Recall from \s\ref{ss:aggp conj} the decomposition of the global Hermitian space $W = W^\flat \oplus Fu$, and recall from \s\ref{sglob int mod} that we have chosen lattices $\Lambda_p \subset W_p$ and $\Lambda_p^\flat \subset W^\flat$ in the respective definitions of $\M_{K_{\wt G}}(\wt G)$ and $\M_{K_{\wt H}}(\wt H)$.  To define the semi-global arithmetic diagonal cycle, we need these lattices to be compatible in a certain way.  If $p$ is split or inert in $F$, or if $p$ is ramified and $n$ is odd, then we assume
\begin{itemize}
\item $\Lambda_p = \Lambda_p^\flat \oplus O_{F_p} u$;
\item if $p$ splits, then $\Lambda_p^\flat$ is self-dual and $\ord_p(u,u) = 0$ (hence $\Lambda_p$ is self-dual);
\item if $p$ is inert, then $\Lambda_p^\flat$ is self-dual and $\ord_p (u,u) = 0$ or $1$ (hence $\Lambda_p$ is respectively self-dual or almost self-dual);
\item if $p$ is ramified and $n$ is odd, then $\Lambda_p^\flat$ is $\pi_p$-modular and $\ord_p (u,u) = 0$ (hence $\Lambda_p$ is almost $\pi_p$-modular); and
\item $K_{\wt H} \subset K_{\wt G} \cap H(\AA_f)$.
\end{itemize}
Note that if $p$ is inert, then we are disallowing the case that $W_p^\flat$ is non-split; see Remark \ref{Wflat nonsplit} below.

As in the case of the Gan--Gross--Prasad Shimura varieties, the inclusion $\wt H \subset \wt G$ is compatible with the Shimura data, and hence induces a morphism of Shimura varieties $\Sh_{K_{\wt H}}(\wt H) \to \Sh_{K_{\wt G}}(\wt G)$.  Still assuming that $n$ is odd if $p$ is ramified, and under the assumptions above, this extends to a morphism of semi-global models
\[
   \M_{K_{\wt H}}(\wt H) \to \M_{K_{\wt G}}(\wt G)
\]
given on points by
\[
   (A_0,A^\flat) \mapsto  (A_0,A^\flat \times A_0),
\]
where the morphism on the additional structure in the moduli problems is similarly defined by ``taking products,'' which we suppress in the notation.  We take the graph morphism
\[
   \Delta\colon \M_{K_{\wt H}}(\wt H) \inj \M_{K_{\wt H}}(\wt H) \times_{\M_0} \M_{K_{\wt G}}(\wt G) = \M_{K_{\wtHG}}(\wtHG),
\]
where the subgroup $K_{\wtHG}$ is as in \eqref{subgroups}.  We define the \emph{semi-global arithmetic diagonal cycle} to be the image in the Chow group,
\[
   z_{K_{\wtHG}} := \Delta_*[\M_{K_{\wt H}}(\wt H)] \in \Ch^{n-1}\bigl(\M_{K_{\wtHG}}(\wtHG)\bigr)_\QQ.
\]
Note that, since we are now working integrally, this Chow group is in the middle codimension.

If $p$ is ramified and $n$ is even, then there is a similar definition of the semi-global arithmetic diagonal cycle, but it is more complicated; see \cite{rsz3}*{\s4.4, \s8.1}.

\subsection{The conjecture}
Now we state the semi-global conjecture.  Keep the assumptions of \s\ref{sglob adc}.  Let $q := \# O_F/(w)$.  Let $\H_{K_{\wtHG}} := C_c^\infty(\wtHG(\AA_f)/\!\!/ K_{\wtHG},\QQ)$ denote the Hecke algebra of smooth bi-$K_{\wtHG}$-invariant \QQ-valued functions with compact support on $\wtHG(\AA_f)$.  Let $f = \otimes_\ell f_\ell \in \H_{K_{\wtHG}}$ be a pure tensor such that, with respect to the product decomposition $\wtHG(\QQ_\ell) = Z(\QQ_\ell) \times (H \times G)(\QQ_\ell)$, $f_\ell$ is of the form $\mathbf{1}_{K_{Z,\ell}^\circ} \otimes \phi_\ell$ for all $\ell$, and $\phi_p = \mathbf{1}_{K_{H,p} \times K_{G,p}}$ (here, as usual, the subscript $\ell$ and $p$ in the subgroup notation means to take the respective $\ell$- or $p$-factor of the subgroup).  Then $f$ acts on the Chow groups of $\M_{K_{\wtHG}}(\wtHG)$ via a Hecke correspondence $R(f)$.  We formally define the intersection number
\begin{equation}\label{int}
   \Int_w(f) := \bigl\la R(f)z_{K_{\wtHG}} , z_{K_{\wtHG}} \bigr\ra \log q,
\end{equation}
where the pairing $\aform$ is defined via the \QQ-linear extension of the Euler--Poincar\'e characteristic $\chi(\,\cdot \otimes^{\mathbf{L}} \cdot\,)$ of the derived tensor product of structure sheaves.  Here is our semi-global arithmetic intersection conjecture \cite{rsz3}*{Conj.\ 8.13}.  We note that part \eqref{afl} is a close variant of an earlier conjecture of W. Zhang \cite{zhang12a} in the context of the GGP Shimura varieties.

\begin{conj}\label{sglob conj}
Let $f$ be as above, and let $f' = \otimes_\ell f_\ell' \in C_c^\infty(G'(\AA))$ be a transfer of $f$ in the sense of \s\ref{m&t}.  Suppose that there is a prime $\ell_0 \neq p$ such that $f$ and $f'$ have regular support at $\ell_0$.
\begin{enumerate}[label=\textup{(\roman*)},ref=\roman*]
\item\label{afl} Assume that $p$ is inert in $F$ and that $\Lambda_p^\flat$ and $\Lambda_p = \Lambda_p^\flat \oplus O_{F_p}u$ are self-dual (and hence $W_p^\flat$ and $W_p$ are split).  Suppose that $f_p' = \mathbf{1}_{G'(\ZZ_{p})}$.  Then
\[
   \Int_w(f) = -\partial J_p(f').
\]
\item\label{atc} In any situation as in \s\ref{sglob adc} (including when $p$ ramifies and $n$ is odd),
\[
   \Int_w (f) = -\partial J_p(f') - J(f'_\corr)
\]
for some correction function of the form $f'_\corr = f_{p,\corr}' \otimes_{\ell \neq p}f_\ell' \in C_c^\infty(G'(\AA_{f}))$.  Furthermore, given $f$, the matching function $f'$ may be chosen such that $f'_\corr = 0$.
\end{enumerate}
\end{conj}

Several remarks are in order.  First, we have already noted that the regular support hypothesis on $f'$ implies that the terms on the right-hand sides of the conjectured equalities are well-defined.  By $f$ having regular support at $\ell_0$, we mean the obvious analog of the condition for $f'$, namely that $\supp f_{\ell_0} \subset \wtHG(\QQ_{\ell_0})_\rs$.  This implies that the intersection of the ``physical'' cycles $|R(f)z_{K_{\wtHG}}| \cap |z_{K_{\wtHG}}|$ is supported in the basic locus of the special fiber \cite{rsz3}*{Th.\ 8.5, Rem.\ 8.6}.  It follows from Theorem \ref{semi-glob model thm} that $\M_{K_{\wtHG}}(\wtHG)$ is a \emph{regular} DM stack in all cases, and therefore the intersection number $\Int_w(f)$ is well-defined.  Furthermore, if $p$ splits, then by loc.\ cit.\ $|R(f)z_{K_{\wtHG}}| \cap |z_{K_{\wtHG}}| = \emptyset$, and by \cite{zhang12a}*{Prop.\ 3.6(ii)} $\partial J_p(f') = 0$.  Thus part \eqref{atc} of the conjecture holds true in the split case and amounts to the equality $0 = 0$ (taking $f_\corr' = 0$), cf.\ \cite{rsz3}*{Th.\ 1.3}.

We further remark that implicit in part \eqref{afl} is that it is possible for $f$ and $f'$ to be transfers under the given hypotheses --- or more precisely, that when $p$ is inert in $F$ and $W_p^\flat$ and $W_p$ are both split, that $\mathbf{1}_{K_Z^\circ,p} \otimes \mathbf{1}_{K_{H_p} \times K_{G,p}}$ transfers to $\mathbf{1}_{G'(\ZZ_p)}$.  This assertion is the (group version of the) \emph{fundamental lemma} conjecture of Jacquet--Rallis \cite{jacral11}, who proposed that the analogous transfer relation holds for any unramified quadratic extension of $p$-adic fields for $p$ odd. In \cite{yun11}, Yun proved the function field version of the fundamental lemma when $p > n$, and Gordon deduced the $p$-adic version from this for $p$ large but unspecified.

There is also an archimedean version of Conjecture \ref{sglob conj} for $p = \infty$, in which $\M_{K_{\wtHG}}(\wtHG)$ is a complex analytic orbifold, using arithmetic Chow groups.  We refer to \cite{rsz3}*{\s8.3} for details.

The evidence we have for the above semi-global conjecture is the following.

\begin{thm}[Zhang \cite{zhang12a}, \cite{rsz3}]
Conjecture \ref{sglob conj} holds for $n \leq 3$.
\end{thm}

Let us briefly comment on the proof.  The case when $p$ splits, which we have already discussed and which is proved for all $n$, is not difficult. But the case of non-split $p$, even for $n \leq 3$, is fundamentally deeper.  There are two main ingredients in the proof.  The first is a kind of local version of the conjecture, where the Shimura variety is replaced by a Rapoport--Zink formal moduli space of $p$-divisible groups: the \emph{arithmetic fundamental lemma} conjecture of Zhang \cite{zhang12a} in the context of part \eqref{afl}, and the \emph{arithmetic transfer} conjecture of \cite{rsz17,rsz18} in the context of part \eqref{atc}.  These have been proved for $n \leq 3$.  The second ingredient is Rapoport--Zink uniformization of the basic locus of the Shimura variety \cite{rapzink96}, which is used to pass from the local statements to the semi-global ones.

\begin{rem}\label{Wflat nonsplit}
Recall that in \s\ref{sglob adc} we disallowed the case that $W_p^\flat$ is non-split when $p$ is inert.  The reason for this is that we do not have an arithmetic transfer conjecture in this case.  We hope to develop such a conjecture (and a proof in low-rank cases) in some future work.
\end{rem}

\section{Global aspects}\label{glob}
To conclude the paper, we sketch the global counterpart to the semi-global conjecture of the previous section; in fact, this serves as motivation for the semi-global conjecture.  We continue to take $F_0 = \QQ$, and we assume that $2$ splits in $F$.

The first step is to define \emph{global} integral models over $\Spec O_F$ of the RSZ Shimura varieties, in the case of appropriate level subgroups.  Let $\Lambda \subset W$ be an $O_F$-lattice whose localization $\Lambda_p$, for each prime $p$, is of the form specified in each case in \s\ref{sglob int mod} according to the type of $(F_p/\QQ_p, W_p)$.\footnote{In particular, when $n$ is even, we are imposing that $W_p$ is split at all ramified primes $p$.}  (It is easy to see that such a $\Lambda$ exists.)  Let
\[
   K_{\wt G}^\circ := K_Z^\circ \times \Stab(\Lambda) \subset \wt G(\AA_f) = Z(\AA_f) \times G(\AA_f).
\]
We define $\M_{K_{\wt G}^\circ}(\wt G)$ to be the moduli problem over $\Spec O_F$ whose points valued in a test scheme $S$ consist of the groupoid of tuples $(A_0,\iota,A,\iota,\lambda)$, where $(A_0,\iota_0) \in \M_0(S)$, $A$ is an abelian scheme over $S$, $\iota\colon O_F \to \End_S(A)$ is an action satisfying the Kottwitz condition \eqref{kott}, and $\lambda\colon A \to A^\vee$ is a polarization, such that for all non-archimedean places $w$ of $F$, the base change of $(A,\iota,\lambda)$ over $S \otimes_{O_F} O_{F,(w)}$ is of the type defined in the semi-global moduli problem for $\wt G$ over $\Spec O_{F,(w)}$ in \s\ref{sglob int mod}.  In particular, this implies that the degree of $\lambda$ is $\#(\Lambda^\vee/\Lambda)$ and $\ker \lambda \subset A[\ddd]$, where $\ddd$ denotes the product of the prime ideals in $O_F$ lying over primes $p$ which ramify in $F$, or which are inert in $F$ and for which $W_p$ is non-split.  In order that the generic fiber of the moduli problem be a single copy of the canonical model of $\Sh_{K_{\wt G}^\circ}(\wt G)$, we also impose the condition that for every geometric point $\ol s \to S$ and prime number $\ell \neq \charac \kappa(\ol s)$, there exists an isomorphism of $O_{F,\ell}/\ZZ_\ell$-Hermitian lattices $\Hom_{O_{F,\ell}}(\rT_\ell(A_{0,\ol s}), \rT_\ell(A_{\ol s})) \simeq -\Lambda_\ell$; see \citelist{\cite{rsz3}*{Rem.\ 5.1(ii)}\cite{bhkry?}*{\s\s2.2--2.3}}.  Here $\rT_\ell(\,\cdot\,)$ denotes the $\ell$-adic Tate module, and the Hom space between Tate modules carries a natural Hermitian pairing as in \eqref{level}.  Then $\M_{K_{\wt G}^\circ}(\wt G)$ is representable by a Deligne--Mumford stack over $\Spec O_F$, whose base change over $\Spec O_{F,(w)}$ is naturally isomorphic to the semi-global integral model for $\Sh_{K_{\wt G}^\circ}(\wt G)$ defined in \s\ref{sglob int mod}, for all non-archimedean places $w$ of $F$ \cite{rsz3}*{Th.\ 5.2}.

Upon choosing an appropriate lattice $\Lambda^\flat \subset W^\flat$, we may analogously define a global integral model $\M_{K_{\wt H}^\circ}(\wt H)$ of $\Sh_{K_{\wt H}^\circ}(\wt H)$.  When $\Lambda^\flat$ and $\Lambda$ are suitably related, and taking toroidal compactifications as before, we obtain a regular flat proper DM stack $\M_{K_{\wtHG}^\circ}(\wtHG)$ over $\Spec O_F$, where $K_{\wtHG}^\circ := K_{\wt H}^\circ \times_{K_Z^\circ} K_{\wt G}^\circ$, together with a global arithmetic diagonal cycle $z_{K_{\wtHG}^\circ} \in \Ch^{n-1}(\M_{K_{\wtHG}^\circ}(\wtHG))$ which is analogous to the semi-global version of \s\ref{sglob adc}.  We remark that the existence of suitable $\Lambda^\flat$ and $\Lambda$ imposes some nontrivial restrictions on the spaces $W^\flat \subset W$, including the dimension constraint that $\lf n/2 \rf$ is odd; see \cite{rsz3}*{Rem.\ 5.3}.  Progress in the direction of Remark \ref{Wflat nonsplit} would allow us to relax some of these restrictions.

Now, for the regular flat proper DM stack $\M_{K_{\wtHG}^\circ}(\wtHG)$, there is defined (unconditionally) the \emph{Gillet--Soul\'e arithmetic intersection pairing} $\sform_\GS$ on the \emph{arithmetic Chow group} $\wh \Ch{}^{n-1}(\M_{K_{\wtHG}^\circ}(\wtHG))$ \cite{gilsoul90,gil09}.  We are going to use this to define intersection numbers, in analogy with \eqref{int}.  We promote $z_{K_{\wtHG}^\circ}$ to a class $\wh z_{K_{\wtHG}^\circ} \in \wh \Ch{}^{n-1}(\M_{K_{\wtHG}^\circ}(\wtHG))$ by adding a Green's current, and we consider a natural Hecke subalgebra $\H_{K_{\wtHG}^\circ}^\spl \subset \H_{K_{\wtHG}^\circ}$ which acts by correspondences $\wh R$ on $\wh \Ch{}^{n-1}(\M_{K_{\wtHG}^\circ}(\wtHG))$ \cite{rsz3}*{\s8.1}.  For  $f \in \H_{K_{\wtHG}^\circ}^\spl$, we define
\[
   \Int(f) := \bigl(\wh R(f)\wh z_{K_{\wtHG}^\circ}, \wh z_{K_{\wtHG}^\circ}\bigr)_\GS.
\]
We then have the following global conjecture \cite{rsz3}*{Conj.\ 8.2}.
\begin{conj}\label{glob conj}
Let $f \in \otimes_\ell f_\ell \in \H_{K_{\wtHG}^\circ}^\spl$ be a pure tensor which transfers to $f' \in C_c^\infty(G'(\AA_{f}))$.  Then
\[
   \Int(f) = -\partial J(f') - J(f'_\corr)
\]
for some correction function $f'_\corr \in C_c^\infty(G'(\AA_{f}))$.  Furthermore, given $f$, the transfer $f'$ may be chosen such that $f'_\corr = 0$.
\end{conj}

As in the introduction, we remark that the distributions on the analytic side of this conjecture are related to the $L$-function in the arithmetic Gan--Gross--Prasad conjecture.

How does the global conjecture relate to the semi-global conjecture?  In fact, aside from a certain formal similarity, it doesn't directly.  Indeed, the global conjecture has the defect that the bi-invariance condition imposed on $f$ by the ``large'' subgroup $K_{\wtHG}^\circ$ is quite restrictive, so that one cannot impose any regular support assumptions on $f$.  In \cite{rsz3}*{\s8.2}, Conjecture \ref{glob conj} is extended to the case of certain deeper level subgroups $K_{\wtHG} \subset K_\wtHG^\circ$, which in terms of the moduli problem correspond to adding Drinfeld level structures at finitely many split primes, cf.\ Remark \ref{drin}.  This allows a larger Hecke algebra to act, and in particular, one may consider transfers $f$ and $f'$ in the conjecture which are regularly supported at some (split, in this case) prime.  For such $f$, the cycles $\wh R(f) \wh z_{K_{\wtHG}}$ and $\wh z_{K_{\wtHG}}$ do not meet in the generic fiber, and the intersection number $\Int(f)$ localizes as a sum $\sum_w \Int_w(f)$ over the nonarchimedean places $w$ of $F$, where $\Int_w(f)$ denotes the semi-global intersection number from before.  Similarly, for such $f'$, we have already noted that there is a localization $\partial J(f') = \sum_p \partial J_p(f')$ in \eqref{delJ decomp}.  In this way, the semi-global conjecture predicts the contributions place-by-place to the global conjecture.

\end{document}